\documentstyle{amsart}
\input{bull-art}
\newcommand{\p}{\partial}

\newtheorem{thm}{Theorem}
\newtheorem{lem}{Lemma}

\newtheorem{conjj}{Conjecture}

\theoremstyle{remark}
\newtheorem{rem}{Remark}
\theoremstyle{remark}

\theoremstyle{remark}
\newtheorem{ex}{Example}

\newcommand{\lr}{\longrightarrow}

\begin{document}
\def\currentvolume{28}
\def\currentissue{1}
\def\currentyear{1993}
\def\currentmonth{January}
\def\copyrightyear{1993}
\def\currentpages{104-108}
\author{Guihua Gong }
\address{Department of Mathematics,
University of Toronto,
Toronto, Ontario,  Canada 
M5S 1A1, and 
Department of Mathematics,
Jilin University, People's Republic of China}
\curraddr{Department of Mathematics, Queen's University, 
Kingston, Canada 
K7L 3N6}
\email{guihua@@math.toronto.edu}
\date{March 24, 1992 and, in revised form, June 25, 1992}
\title {Relative $K$-cycles and Elliptic Boundary 
Conditions}
\dedicatory{Dedicated to Professor Zejian Jiang on his 
seventieth birthday}

\maketitle

\begin{abstract}
In this paper, we discuss the following conjecture raised 
by Baum-Douglas:  For
any first-order elliptic differential operator $D$ on 
smooth manifold $M$ with
boundary $\p M$, $D$ possesses an elliptic boundary 
condition if and only if
$\partial [D]$ = 0 in $K_1(\partial M)$, where $[D]$ is 
the relative $K$-cycle 
in $K_0(M, \partial M)$ corresponding to $D$.  We prove 
the ``if'' part of this
conjecture for $\dim(M)$ $\not=$ 4, 5, 6, 7 and the ``only 
if'' part of the
conjecture for arbitrary dimension.   
\end{abstract}

First we  fix some notation.  $M$ is a compact oriented 
smooth manifold with 
smooth boundary $\partial M$.  We always suppose that $M$ 
is embedded in some
compact smooth manifold $\widetilde M$ without boundary  
of the same dimension 
(e.g., $\widetilde M$ can be taken as double of $M)$.  We 
denote $\stackrel{
\circ}{M} = M \setminus \partial M$.   Furthermore, we 
assume that $E_0$ and
$E_1$ (in fact, all the vector bundles in this paper) are 
smooth complex
Hermitian vector bundles over $M$ and  that  $D: 
C^{\infty}(E_0) \rightarrow
C^{\infty}(E_1)$ is a first-order elliptic differential 
operator from smooth
sections of $E_0$ to that of $E_1$.   By $H^s(M, E_i)$ and 
$H^s(\partial M,
E_i)$ we shall denote the Sobolev spaces of sections of 
$E_i$ and
$E_i|_{\partial M}$ with respect to fixed smooth measures 
on $M$ and $\partial
M$, respectively.  

The elliptic boundary value problem (an elliptic operator 
with an elliptic
boundary condition) has been studied for a long time.  As 
noted in
[1, 5, 6] and other references, there exist 
topological obstructions to  impose an elliptic boundary 
condition on the above
$D$.  A fundamental problem is to find all such 
obstructions.  
Baum, Douglas, and 
Taylor constructed a relative K-cycle  $[D] \in K_0(M, 
\partial M)
\cong KK(C_0(\stackrel{\circ}{M}),{\Bbb C})$ (here  
$C_0(\stackrel{
\circ}{M})$ is the algebra of continuous functions on $M$ 
which vanish on
$\partial M)$ corresponding to $D$ (see [2--4]
for details). From the definition of  relative K-homology 
group $K_0(M,
\partial M)$ given by Baum, Douglas, and Taylor, the 
boundary map  
$\partial\!: K_0(M,
\partial M) \longrightarrow K_1(\partial M)$ is very 
concrete 
[2--4].
Also Baum and Douglas conjectured that the only
obstruction for $D$ possessing elliptic boundary  
conditions is that $\p [D]
\neq 0$.  More precisely, the following conjecture  first 
appeared in
[2] in a closely related form.  

\begin{conjj}
There exist a vector bundle $E_2$ over  $\partial M$ and a
zeroth-order pseudo-differential operator $B$ defined from 
$C^{\infty}(\partial
M, E_0)$ to $C^{\infty}(\partial M, E_2)$ such that 
$$
\begin{array}{l}
 \left( \begin{array}{c}
          D\\
          B \circ \gamma \end{array}  \right):\ \   H^1(M, 
E_0) \longrightarrow 
\begin{array}{c}
H^0(M, E_1)\\
\oplus\\
H^{1/2}(\partial M, E_2)
\end{array}
\end{array}
$$
is Fredholm if and only if $\p [D] = 0$ in $K_1(\partial 
M)$,
where $\gamma:  H^1(M, E_0) \longrightarrow 
H^{1/2}(\partial 
M, E_0)$ is the trace map.  
\end{conjj}

\begin{rem}
Let $D$ be as above with pincipal symbol $p(x, \xi)$.  A 
zeroth-order
pseudo-differential operator $B$ with principal symbol 
$b(x, \xi)$ from
$C^{\infty}(\p M, E_0)$ to  $C^{\infty}(\p M, E_2)$ is 
said to be elliptic to
$D$ (see [6], p.\ 233) if, for every $x \in \p M$ and $\xi 
\in
T_x^*(\p M)$, the map $M_{x,\xi}^+ \ni u \lr b(x, \xi)u(0) 
\in (E_2)_x$ is
bijective, where $T_x^*(\p M)$ and $(E_2)_x$ are the 
fibres at the point $x$ of
the cotangent bundle $T^*\p M$ and the bundle $E_2$,  and, 
furthermore, $M_{x,
\xi}^+$ is the set of all $u \in C^{\infty}({\Bbb 
R},(E_0)_t)$ with $p(x, \xi -i
\frac{d}{dt}\cdot n_x)u(t) = 0$ $(n_x$ is the interior 
conormal vector  of $M$
at $x)$ which are bounded on ${\Bbb R}_+$.  If $B$ is 
elliptic to $D$, then the
above 
\begin{math} \bigl( \begin{smallmatrix}
          D\\
          B \circ \gamma 
\end{smallmatrix} \bigr) 
\end{math}
is Fredholm.  Such a system 
\begin{math} \bigl( \begin{smallmatrix}
          D\\
          B \circ \gamma 
\end{smallmatrix} \bigr) 
\end{math}
is often called an elliptic boundary value problem or an 
elliptic operator with
an elliptic  boundary condition; meanwhile, $D$ is also 
said to possess an
elliptic boundary condition.
\end{rem}

\begin{rem}
Although the above elliptic boundary condition is used in 
most references, the
original form of the conjecture in [2] is in a slightly 
different form
from the above. In [2] the operator for the boundary 
condition is of the
form $\gamma \circ B$, where $B$ is a zeroth-order 
pseudo-differential operator
from $E_0$ to a smooth vector bundle over a neighborhood 
of $M$. The reason we
use a slightly different form of the conjecture is as 
follows:  for general
zeroth-order pseudo-differential operator $B$ defined on 
$\widetilde M$, there
is no canonical way to restrict $B$ to $M$ as an operator  
$B_{M}: H^s(M)
\longrightarrow H^s(M)$ when $s > 0$.  So one needs to put 
some restriction on
$B$.  One of the natural restrictions  is that $B$ has the 
transmission
property with respect to $\p M$ (see [5]).   We also prove 
our theorem
for this kind of boundary condition (see Theorem 1).  It 
must be pointed
out that the existence of a boundary condition of type $B 
\circ \gamma$ implies
the existence of that of type $\gamma \circ B$. 
\end{rem}

In this paper, we prove the ``only if'' part of the 
conjecture which can be
thought of as a generalization of Corollary 4.2 in [4] 
(there $B$ is a
differential operator).  Conversely, we prove that if 
$\dim(M) \neq 4, 5, 6, 7$
and $\p [D] = 0$, then $D$ possesses an elliptic boundary 
condition as in
Remark 1. Hence the ``if'' part of the conjecture has been 
proved for $\dim(M)
\neq 4, 5, 6, 7.$  The  cases of $\dim(M)$ being equal to 
$4, 5, 6$,
 or 7 are still
open, but we prove a theorem which can be thought of as 
the ``if'' part of the
conjecture in the sense of stablization in K-homology 
group for arbitrary
dimension. Our results will be useful for constructing 
absolute $K$-cycles in
$K_0(M)$ which are preimages of $[D] \in K_0(M, \p M)$ 
under the canonical map
from $K_0(M)$ to $K_0(M, \p M)$ when $\p [D] = 0$.

Our main results are the following:

\begin{thm} \RM{(``}only if\ \RM{''} part\RM)  
$\p [D] = 0$  if one of the following is true\RM:

\RM{(i)}  There exist a smooth vector bundle $E_2$ over 
$\p M$ and a 
zeroth-order pseudo-differential operator $B$ from 
$E_0|_{\p M}$ to $E_2$ 
such that 
\begin{math} \bigl( \begin{smallmatrix}
          D\\
          B \circ \gamma 
\end{smallmatrix} \bigr) 
\end{math}
in the conjecture is Fredholm. 

\RM{(ii)} There exist a bundle $E_2$ over a neighborhood 
of $M$ in ${\widetilde
M}$ and a zeroth-order pseudo-differential operator$B$ 
with transmission
property with respect to $\p M$ from $E_0$ to $E_2$ such 
that 
$$
\begin{array}{l}
 \left( \begin{array}{c}
          D\\
          \gamma \circ B \end{array}  \right):\ \   H^1(M, 
E_0) \longrightarrow 
\begin{array}{c}
H^0(M, E_1)\\
\oplus\\
H^{{1}/{2}}(\partial M, E_2)
\end{array}
\end{array}
$$
is Fredholm.
\end{thm}

\begin{thm}
\RM{(``}if\ \RM{''} part\RM)  If $\p [D] = 0$, then there 
exists a first-order
elliptic differential operator $D_1$ acting on smooth 
vector bundles over $M$
with $[D_1] = 0$ in $K_0(M, \p M)$ such that $D\oplus D_1 
$ possesses an
elliptic boundary condition as in Remark \RM1, and, 
furthermore, if $\dim (M) \neq
4, 5, 6, 7$, then $D$ itself possesses an elliptic  
boundary condition.  
\end{thm}

The main idea of the proof of Theorem 1 is to construct an 
intertwining between
$\p [D]$ and a trivial element in $K_1(\p M)$.  In the 
proof, we use Calderon
projection, functional calculus of pseudo-differential 
operators (including
Boutet de Monvel type operators), and the techniques in 
the proof of Proposition
4.5 of [4].

The proof of Theorem 2 makes use of two key lemmas (see 
below).

Let $ST^*\p M$ be the unit sphere bundle of $T^*\p M$ over 
$\p M$ and $\pi\!:
ST^*\p M \lr \p M$ be the canonical projection map.  Let 
$\widetilde E_0 = \pi^
*(E_0|_{\p M})$ be the bundle over $ST^*\p M$.  We write 
the principal symbol
of $D$, in a coordinate neighborhood $U$ of $x \in \p M$, 
as 
$$p(x, x_n, \xi, \xi_n) = \sum_{j=1}^{n-1}p_j(x, x_n)\xi_j 
+ p_n(x,x_n)\xi_n,$$ 
where  $x_n$ is the coordinate for the normal direction of 
$\p M$. We define
$$\tau(x, \xi) = ip_n^{-1}(x, 0)\sum_{j=1}^{n-1}p_j(x, 
0)\xi_j$$
for $x \in \p M$ and $\xi \in ST^*\p M$.  Then $\tau(x, 
\xi)$ is a map from a
fibre of $E_0$ into itself and has no purely imaginary 
eigenvalue. Let
$V_{\pm}$ be the subbundle of ${\widetilde E_0}$ over 
$ST^*\p M$ corresponding
to the span of the generalized eigenvectors of $\tau(x, 
\xi)$ corresponding to
the eigenvalues with positive/negative real parts.

\begin{lem}
$\p [D] = 0$ if and only if $[V_+] \in \pi^*K^0(\p M) 
\subset K^0(ST^*\p M)$.
\end{lem}

\begin{lem} If $E_0 $ and $E_1$ are vector bundles over 
$M$ which allow a
first-order elliptic differential operator $D$ to act from 
one to the other,
and if $\dim(M) = n$, then 

\RM{(i)} $f\dim(E_0) = f\dim E_1 \geq 2^{[({n-1})/{2}]}$\RM;

\RM{(ii)} $f\dim(E_0) = f\dim E_1 \geq 2^{[({n-1})/{2}]+
1}$ provided $n$ is
even and $\p [D] = 0 $,

where $f\dim$ denotes dimension of each fibre of the 
vector bundles.
\end{lem}

\begin{proof*}{The Proof of Theorem \RM 2}
By Lemma 1, if $\p [D] = 0 $, one has 
\[[V_+] \in \pi^*K^0(\p M) \subset
K^0(ST^*\p M).\]
 By Lemma  2, $f\dim (V_+) = \frac{f\dim (E_0)}{2} \geq
2^{[{n}/{2}]-1}$. Therefore, $f\dim (V_+) \geq  n-1 > 
\frac{\dim (ST^*\p
M)}{2} = \frac{2n-3}{2}$, whenever  $\dim (M) \geq 8$. 
Hence there exists a
complex vector bundle $E_2$ over $\p M$, such that $V_+ 
\cong \pi^*E_2$. This
is also true for $\dim (M) \leq 3$, since the collection 
of complex vector
bundles over $ST^*\p M$  
$(\dim (ST^*\p M) \leq 3)$  has property of
cancellation.

Let $\psi$ be the bundle isomorphism
$$
\begin{array}{ccc}
V_+        & \buildrel \psi \over \lr & \pi^*E_2\\
\downarrow &                         &   \downarrow \\
ST^* \p M  & \lr                     & ST^* \p M
\end{array}
$$

For any $(x, \xi) \in ST^* \p M$, let $b(x, \xi)$ be the 
bundle map defined by
$$
\begin{array}{ccccc}
E_0      &  \stackrel{\text{project to}}{\lr}  &  V_+  
&  \stackrel{\psi}{\lr}   &    E_2 \\
\downarrow  &  &   &   &   \downarrow \\
ST^* \p M &   \lr  &   ST^* \p M   &    \lr  &  ST^* \p M 
\end{array}
$$

Furthermore, let $B :  E_0|_{\p M}  \lr  E_2$  be the 
zeroth-order
pseudo-differential operator with symbol $b(x, \xi)$.   It 
follows that  $B$ is
elliptic to $D$ (see [6]). \end{proof*}

\begin{ex} 
If $M$ is a $spin^c$ manifold with smooth boundary and $D$ 
is the 
Dirac operator over $M$, then it  is computed in [4] that 
$\p[D] \neq
0$. Hence $D$ possesses no elliptic boundary condition 
(even possesses no
boundary condition as in the conjecture).
\end{ex}

\begin{ex}
For any $D$, let $D^*$ be the formal adjoint of $D$.  It is
easy to prove that $[D] = -[D^*]$ in $K_0(M,\p M)$; hence, 
$\p [D\oplus D^*]
=0$. It follows that $D\oplus D^*$ possesses an elliptic 
boundary condition
provided $\dim (M) \neq$ 4  or 6.  (It should be noted 
that we only need to
exclude the manifolds with dimension 4 and 6 here, since 
the dimension of
the bundle on which  $D\oplus D^*$ acts is twice the 
dimension of the bundle on
which $D$ acts.)
\end{ex}

The details of the proofs will appear elsewhere.

\section*{Acknowledgments} This work was done while the 
author was a
postdoctoral fellow at the University of Toronto.  He is 
grateful to Professors
George A. Elliott and  Man-Deun Choi for their support. 
The author also thanks
Professors Man-Deun Choi, Ronald G. Douglas, George A. 
Elliott,  and Peter
Greiner and Ms. Liangqing Li for many helpful conversations.


\begin{thebibliography}{6}

\bibitem[1]{APS} M. Atiyah, V. Patodi, and 
I. Singer, {\em Spectral asymmetry and
riemannian geometry}. I, II, III, 
 Math. Proc. Cambridge Philos. Soc.  {\bf 77} (1975), 
43--69;
{\bf 78} (1975), 405--632; {\bf 79} (1976),  71--99.
\bibitem[2]{BD1} P. Baum  and R. Douglas, {\em Index 
theory, bordism, and
K-homology,} Operator Algebras and $K$-Theory (R. G. 
Douglas and C. Schochet,
eds.), 
Contemp. Math., vol. 10, Amer. Math. Soc., Providence, RI, 
1982,
pp.  1--31.
\bibitem[3]{BD2} \bysame,  {\em Relative $K$-homology and
$C^*$-algebra,} manuscript.

\bibitem[4]{BDT} P. Baum, R. Douglas, and M. Taylor, {\em 
Cycles and relative
cycles in analytic $K$-homology,} J. Differential Geom. 
{\bf 30} 
(1989), 761--804.

\bibitem[5]{Bout} L. Boutet de Monvel, {\em Boundary 
problems for
pseudodifferential operators,} Acta Math. {\bf 126} 
(1971), 11--51.

\bibitem[6]{Hom} L. H\"omander, {\em The analysis of 
linear partial
differential operators}. III, Springer, New York, 1985.

\end{thebibliography}
\end{document}